\documentclass{article}

\usepackage[english]{babel}
\usepackage{appendix}

\usepackage[letterpaper,top=2cm,bottom=2cm,left=3cm,right=3cm,marginparwidth=1.75cm]{geometry}

\usepackage{setspace}
\usepackage{amsmath}
\usepackage{amssymb}
\usepackage{graphicx}
\usepackage[colorlinks=true, allcolors=blue]{hyperref}
\usepackage{authblk}
\numberwithin{equation}{section}
\newtheorem{theorem}{Theorem}[section]
\newtheorem{lemma}[theorem]{Lemma}
\newtheorem{proposition}[theorem]{Proposition}
\newtheorem{corollary}[theorem]{Corollary}

\newcommand{\scal}{\mathrm{scal}}
\newcommand{\Ric}{\mathrm{Ric}}

\title{On Shrinking Ricci solitons with positive isotropic curvature in higher dimensions}
\author[1]{Zhengnan Chen}
\affil[1]{School of Mathematical Sciences, Peking University, Beijing 100871, China}
\date{}
\begin{document}
\maketitle
\footnotetext[1]{E-mail address: \href{mailto:znchen@pku.edu.cn}{znchen@pku.edu.cn}}
\begin{abstract}

For all dimensions $n\geq5$, let $(M,g,f)$ be a $n-$dimensional shrinking gradient Ricci soliton with strictly positive isotropic curvature (PIC). Suppose furthermore that $\nabla^2f$ is $2-$nonnegative and the curvature tensor is WPIC1 at some point $\bar{x}\in M$. Then $(M,g)$ must be a quotient of either $S^{n}$ or $S^{n-1}\times\mathbb{R}$.  Our result partially extends the classification result for 4-dimensional PIC shrinking Ricci solitons established in \cite{LNW16} to highter dimensions. Combining the pinching estimates deduced in \cite{chen25}, we also extend the result in \cite{cho2023} to dimensions $n\geq9$. Namely that a complete ancient solution to the Ricci flow of dimension $n\geq9$ with uniformly PIC must be weakly PIC2.
\end{abstract}


\section{Introduction}
A shrinking gradient Ricci soliton (SGRS for short) is a triple $(M,g,f)$, consisting of a complete Riemannian manifold $(M,g)$ and a smooth function $f$ (called the potential function) satisfying the following equation
\begin{equation}\label{SGRS eq}
    \Ric+\nabla^2f=\frac{1}{2}g.
\end{equation}

Under Ricci flow, a SGRS evolves in a self-similar and shrinking way. From equation \ref{SGRS eq} we see that it generalize Einstein manifolds with positive curvature. The study of SGRS plays an important role in understanding the formation of singularities in the Ricci flow.

There are many works regarding classification of SGRSs in low dimensions. In dimension 2, by works of \cite{Hamilton88}, any nonflat SGRS must be a shrinking sphere $S^2$ or its quotient $\mathbb{R}P^2$. In dimension 3, Perelman \cite{Perelman02} firstly classified $\kappa$-noncollapsed SGRS with bounded nonnegative sectional curvature. Later classification works \cite{Naber10}\cite{NW08} and \cite{CCZ08} gradually removed the non-collapsing and all curvature assumptions. By their works, any 3 dimensional SGRS must be a quotient of either $S^3$ or $S^2\times\mathbb{R}$.

In dimension 4 and higher dimensions, it is difficult to classify all SGRSs. To the author's knowledge, existing classification results to date all require certain assumptions on the curvature tensor. In \cite{MW15}, Munteanu and Wang showed that a SGRS with nonnegative sectional curvature and positive Ricci curvature must be compact. Combining this result with the pinching cone estimate established in \cite{BohmWilking08}, we can completely classify SGRSs with nonnegative curvature operator. In dimension 4 there are various results assuming different curvature assumptions in recent years. Interested readers can see section 4.3 of \cite{cao2024geometry} for reference.

In this paper, we are mainly interested in the positive isotropic curvature (PIC for short, see section 2.2 for relevant definition) condition, which was introduced by Micallef and Moore in \cite{MicallefMoore88} in their study of minimal two-spheres in Riemannian manifolds. This curvature condition makes sense in dimension $n\geq4$ and plays an important role in the theory of Ricci flow. In \cite{Hamilton97FourmanifoldsWP}, Hamilton initiated and developed the study of the Ricci flow starting from PIC manifolds (manifolds satisfying PIC condition) in dimension 4. He showed that PIC is preserved along Ricci flow in dimension 4, and obtained pinching estimate for 4 dimensional PIC condition. Later in \cite{BrendleSchoen07}, Brendle and Schoen showed that PIC and its variants PIC1, PIC2 are all preserved along Ricci flow in all dimensions $n\geq 4$ (note that \cite{Nguyen10} also independently showed that PIC is preserved along Ricci flow in all dimensions).  

Using estimates in \cite{Hamilton97FourmanifoldsWP} and applying an integration by parts argument, \cite{LNW16} classified 4 dimensional strictly PIC SGRSs. Such solitons must be quotient of either $S^4$ or $S^3\times \mathbb{R}$. See \cite{ni20074}, \cite{richard2016positive}, \cite{cao2023four} for related works.

In higher dimensions, \cite{Li2019} gave a classification for weakly PIC1 SGRSs. For the PIC case, in dimension $n\geq12$, \cite{Brendle19} used the pinching cone method (such method was primarily introduced in \cite{BohmWilking08}) to construct a family of pinching cones from PIC towards weakly PIC1. Recently, the author in \cite{chen25} improved the above pinching cone estimate to dimension $n\geq 9$. These pinching estimates have important topological applications regarding topology of PIC manifolds, see \cite{Brendle19},\cite{chen25},\cite{huang2019} and \cite{huang2023}. Keaton Naff in \cite{naff19} used the family of pinching cones constructed in \cite{Brendle19} to show that all nonflat uniformly PIC SGRS with dimension $n\geq12$ must be weakly PIC2. More generally, in the ancient solution setting, Jae Ho Cho and Yu Li showed in \cite{cho2023} that a nonflat noncompact complete ancient solution of dimension $n\geq12$ that is uniformly PIC must be weakly PIC2.

We can use the pinching estimate deduced in \cite{chen25} to show that the above result holds actually for dimension $n\geq9$.

\begin{theorem}\label{thm11}
    Let $(M,g(t)),t\in(-\infty,0]$ be a complete ancient solution to the Ricci flow of dimension $n\geq9$ with uniformly PIC. Then it is weakly PIC2.
\end{theorem}

The proof follows directly form the proof in \cite{cho2023}. When we run the arguments in \cite{cho2023}, the only difference is that the family of pinching cones we use is different. Therefore it suffices to show that the (0,4)-curvature tensor $R$ of $(M,g)$ is contained in $\hat{C}(b)$ for some $b>0$ (the pinching family $\{\hat{C}(b)\}$ is stated in Theorem 1.2 of \cite{chen25}). We will introduce the concept of family of pinching cones and briefly show Theorem \ref{thm11} at the end of section 2.

Applying Theorem 1.1 of \cite{naff19} in SGRS case and Theorem 5.7 of \cite{cho2023} in ancient solution case, we obtain the below corollary:

\begin{corollary}
    \begin{enumerate}
        \item Let $(M,g,f)$ be a non-flat SGRS of dimension $n\geq9$ with uniformly PIC. Then $(M,g)$ is isometric to quotient of either $S^n$ or $S^{n-1}\times\mathbb{R}$.
        \item Let $(M,g(t)),t\in(-\infty,0]$ be a complete, $\kappa-$noncollapsed, noncompact ancient solution to the Ricci flow with uniformly PIC. Then $(M,g(t))$ is isometric to the Bryant soliton or a quotient of the shrinking cylinder $S^{n-1}\times\mathbb{R}$.
    \end{enumerate}
\end{corollary}

Given the classification result above and the result of \cite{LNW16} in dimension 4, it is natural to ask can we classify strictly PIC SGRSs, instead of uniformly PIC SGRSs (this is also questioned in \cite{naff19}) for all dimension $n\geq5$. Here we give a partial answer:  

\begin{theorem}\label{main thm}
    Suppose $n\geq5$, $(M,g,f)$ is a non-flat $n-$dimensional SGRS satisfying equation \ref{SGRS eq} with curvature operator that is strictly PIC. Assume further that $\nabla^2f$ is $2-$nonnegative (i.e. the sum of 2 smallest eigenvalues of $\nabla^2f$ is nonnegative), and that the curvature tensor $R$ is WPIC1 at some point $\bar{x}\in M$.
    
    Then $(M,g,f)$ must be a quotient of $S^n$ or $S^{n-1}\times\mathbb{R}$.
\end{theorem}

Since now we are not assuming uniformly PIC, we are unable to directly show that the curvature tensor is contained in $\hat{C}(b)$ for some $b>0$. The method we use here is to construct a nonpositive function $p_1(x)$ that measures the deviation of the curvature tensor form the cone $C_{PIC1}$ (see the proof of Proposition \ref{propPICPIC1}  ). Then we deduce a Laplacian type differential inequality for $p_1$. Using a maximal principle argument, we show that under the assumptions above, $p_1\equiv0$, which means $(M,g)$ is WPIC1. Now by Proposition 6.2 in \cite{Li2019}, we conclude that $(M,g)$ is actually WPIC2. The assumption that $\nabla^2f$ is $2-$nonnegative, and that the curvature tensor at at least one point is WPIC1 are satisfied by both $S^n$ and $S^{n-1}\times\mathbb{R}$. These assumptions are used in our analysis of the function $p_1$. It's interesting to ask whether these additional assumptions can be removed.

The rest of this article is organized as follows: In section 2, we recall some standard results for SGRSs and positive isotropic curvature. We prove Theorem \ref{thm11} at the end of section 2. In section 3, we give proof to Theorem \ref{main thm}.

\textbf{Acknowledgements.} The author is grateful to his advisor professor Gang Tian for his constant encouragements. He also thanks Professor Shijin Zhang for valuable discussions.

\section{Preliminaries and Proof of Theorem \ref{thm11}}

\subsection{Standard results for SGRSs}
A SGRS (Shrinking gradient Ricci soliton) would satisfy certain standard equations, which are deduced from the equation \ref{SGRS eq}. One can see for example \cite{Chow2023} for reference.
\begin{equation*}
    \begin{split}
        \scal+\Delta f&=\frac{n}{2},\\
        \Ric(\nabla f,\cdot)&=\mathrm{div}(\Ric)=\frac{1}{2}\nabla\scal,\\
        \mathrm{d}(\scal+|\nabla f|^2-f)&=0.
    \end{split}
\end{equation*}

After shifting $f$ by a constant if necessary, we may assume that 
\begin{equation*}
    \scal +|\nabla f|^2=f.
\end{equation*}

A SGRS satisfying the assumptions above is called normalized. We have the following elliptic equation for normalized SGRSs:

\begin{equation*}
    \begin{split}
        \Delta_f R&= R-Q(R),\\
        \Delta_f\Ric_{ij}&=\Ric_{ij}-2\sum_{k,l}R_{ikjl}\Ric_{kl},\\
        \Delta_f \scal&=\scal-2|\Ric|^2.
    \end{split}
\end{equation*}

Here, the $f-$Laplacian is defined by $\Delta_f=\Delta-\nabla_{\nabla f}$. The quadratic $Q(R)$ is given by $Q(R)=R^2+R^{\#}$, where
    \begin{equation*}
    \begin{split}
        (R^2)_{ijkl}&=\sum_{p,q}R_{ijpq}R_{klpq},\\
        (R^{\#})_{ijkl}&=2\sum_{p,q}(R_{ipkq}R_{jplq}-R_{iplq}R_{jpkq}).
    \end{split}
    \end{equation*}

Notice that standard calculations show that the curvature tensor $R$ evolves by $\frac{\partial R}{\partial t
}=\Delta R+Q(R)$ under evolving frames (See \cite{MorganTian07} Chapter 3 for reference). The elliptic equations above are deduced from the evolution equation for the associated Ricci flow of $(M,g,f)$.

\subsection{Positive isotropic curvature and pinching cones}

We follow the notations as in \cite{chen25}. An algebraic curvature operator $R$ is said to be WPIC (weakly positive isotropic curvature ) if for any orthonormal four-frame $\{e_1,e_2,e_3,e_4\}$, we have
\begin{equation}
    R_{1313}+R_{1414}+R_{2323}+R_{2424}-2R_{1234}\geq 0.
\end{equation}
$R$ is said to be WPIC1 if for any orthonormal four-frame $\{e_1,e_2,e_3,e_4\}$ and any $\lambda\in [0,1]$,
\begin{equation}
    R_{1313}+\lambda^2R_{1414}+R_{2323}+\lambda^2R_{2424}-2\lambda R_{1234}\geq 0.
\end{equation}
Moreover $R$ is said to be WPIC2 if for any orthonormal 4-frame $\{e_1,e_2,e_3,e_4\}$ and any $\lambda,\mu\in [0,1]$.
\begin{equation}
    R_{1313}+\lambda^2R_{1414}+\mu^2R_{2323}+\lambda^2\mu^2 R_{2424}-2\lambda\mu R_{1234}\geq 0.
\end{equation} 

We call $C_{PIC}, C_{PIC1}, C_{PIC2}$ the cones consisting of algebraic curvature operators that satisfy WPIC (weakly PIC), WPIC1 or WPIC2 respectively. Moreover, $R$ is said to be PIC, PIC1 or PIC2 if it's contained in the interior of $C_{PIC1},C_{PIC1}$ or $C_{PIC2}$ respectively. A Riemannian manifold $(M,g)$ is a (strictly) PIC, PIC1 or PIC2 manifold if its associated curvature operator is PIC, PIC1 or PIC2 at every point $x\in M$ respectively.  

Let $S_B^2(\mathfrak{so}(n))$ denote the set of algebraic curvature operators. Let $J\subset\mathbb{R}$ be an interval. A collection of conic sets $\{\{C(b)\}_{b\in J}\}$ parametrized by $b\in J$ is called a family of pinching cones if for each $b\in J$, $C(b)$ satisfies the following conditions:

\begin{enumerate}
    \item $C(b)$ is a closed, convex, $O(n)-$invariant conic subset of $S_B^2(\mathfrak{so}(n))$;
    \item Any $R\in C(b)\backslash\{0\}$ has positive scalar curvature;
    \item  $C(b)$ is transversally invariant under the Hamiltonian ODE $\frac{d}{dt}R=Q(R)$, which means that for any $R\in\partial C(b)\backslash\{0\}$, the quadratic $Q(R)$ lies in the interior of the tangent cone $T_{R}(C(b))$.
    \item The identity $I$ lies in the interior of $C(b)$, where $I=\frac{1}{2}\mathrm{id}\wedge \mathrm{id}$, here the symbol $\wedge$ denotes the Kulkarni-Nomizu product:
    \begin{equation*}
        (A\wedge B)_{ijkl}=A_{ik}B_{jl}+A_{jl}B_{ik}-A_{il}B_{jk}-A_{jk}B_{il}.
    \end{equation*}
\end{enumerate}

Now we let $J=[0,1]$ be a closed interval. And $\{C(b)\}_{b\in[0,1]}$ be a family of pinching cones.
In Theorem 3.3 of \cite{cho2023}, Jae Ho Cho and Yu Li showed that for a complete ancient solution $(M^n,g(t))_{t\in(-\infty,0]}$, if the curvature tensor $R(x,t)$ lies in $C(0)$ for all $(x,t)\in M\times(-\infty,0]$, then $R(x,t)$ lies in $C(1)$ for all $(x,t)\in M$. 

Using the family of pinching cones constructed in \cite{Brendle19}, they showed that in dimension $n\geq12$, any complete ancient solution to Ricci flow with uniformly PIC must be WPIC2. Combining the author's recent result in \cite{chen25}, we can actually show this for $n\geq9$. It suffices to show that the curvature tensor $R$ lies in $C(b)$ for some small $b>0$, where $C(b)$ is defined as in Definition 3.1 of \cite{chen25}.

\subsection{Proof of Theorem \ref{thm11}}

Let $n\geq9$ and $(M^n,g)$ be a complete ancient solution to Ricci flow with uniformly PIC, namely $R-\theta\scal I\in C_{PIC}$ for some $\theta>0$. By Proposition 3.1 of \cite{cho2023}, there exist a uniform $\delta=\delta(n,\theta)>0$ such that
\begin{equation*}
    \lambda_1+\lambda_2\geq\delta \scal
\end{equation*}
holds for all $(x,t)\in M\times(-\infty,0]$. Here $\lambda_1+\lambda_2$ is the sum of two smallest eigenvalues of $\Ric$. 

We verify that for some small $b>0$, $S=l_{a,b}^{-1}(R)$ lies in $\mathcal{E}(b)$ stated in Definition 3.1 of \cite{chen25}. Here $a=\frac{(2+(n-2)b)^2}{2(2+(n-3)b)}b$ goes to $0$ as $b$ goes to $0$. The linear transformation $l_{a,b}$ was introduced in \cite{BohmWilking08}. It acts on space of algebraic curvature operators by

\begin{equation*}
    l_{a,b}(R)=R+b\ \Ric\wedge\mathrm{id}+\frac{2}{n}(a-b)\scal\ I.
\end{equation*}

Hence $l_{a,b}$ varies continuously on $b$. By continuity, we can choose $b$ small enough (which depends on $\theta$ and $n$), such that $S-\frac{\theta}{2}\scal(S)I\in C_{PIC}$ and $\Ric(S)_{11}+\Ric(S)_{22}\geq0$. Now choose $T=\frac{\theta}{2}\scal(S)I$, then $T$ is positive definite and $S-T\in C_{PIC}$. Notice that the data $\omega$ defined in the beginning of section 3 of \cite{chen25} satisfies $\omega\rightarrow+\infty$ as $b\rightarrow0$. Hence by choosing a smaller $b$ if necessary, we can let the below hold

\begin{equation*}
    \Ric(S)_{22}-\Ric(S)_{11}\leq \omega^{\frac{1}{2}}(\scal(S))^{\frac{1}{2}}\left(\sum_{p=3}^n(T_{1p1p}+T_{2p2p})\right)^{\frac{1}{2}}
\end{equation*}
for every orthonormal frame $\{e_1,\cdots,e_n\}$. Hence, the curvature tensor $S$ lies in $\mathcal{E}(b)$, and thus $R=l_{a,b}(S)$ lies in $C(b)$. This completes the proof of Theorem \ref{thm11}.

\section{Proof of Theorem \ref{main thm}}

We will first need a special case of Proposition A.8 of \cite{Brendle19}.

\begin{lemma}\label{Lem for Q(R)PIC1}
    Let $R\in S_B^2(\mathfrak{so}(n))$ be an algebraic curvature tensor, $c$ be a fixed real, suppose that 
    \begin{equation*}
        Z=R_{1313}+\lambda^2R_{1414}+R_{2323}+\lambda^2R_{2424}-2\lambda R_{1234}+c(1-\lambda^2)\geq0
    \end{equation*}
for every orthonormal four-frame $\{e_1,e_2,e_3,e_4\}$ and every $\lambda\in [0,1]$. Moreover, suppose that $Z=0$ for one particular orthonormal four-frame $\{e_1,e_2,e_3,e_4\}$ and one particular $\lambda\in [0,1)$, then

\begin{equation*}
    Q(R)_{1313}+\lambda^2Q(R)_{1414}+Q(R)_{2323}+\lambda^2Q(R)_{2424}-2\lambda Q(R)_{1234}+c(1-\lambda^2)(\Ric_{11}+\Ric_{22})\geq0.
\end{equation*}
     
\end{lemma}

\textit{Proof.} This is done by setting $H=\frac{c}{2}\mathrm{id}$ in Proposition A.8 of \cite{Brendle19}. We present the proof here for completeness.

Define a new tensor $T\in S_B^2(\mathfrak{so}(n+1))$ by
\begin{equation*}
    T_{ijkl}=R_{ijkl},\quad T_{ijk0}=0,\quad T_{i0k0}=\frac{c}{2}\delta_{ik},
\end{equation*}
for $i,j,k,l\in \{1,2,\cdots,n\}$. It's easily checked that $T$ is well-defined and satisfies the first Bianchi identity. We claim that $T\in C_{PIC}$. That is, for any linearly independent $\tilde{\zeta},\tilde{\eta}\in \mathbb{C}^{n+1}$ that satisfy $g(\tilde{\zeta},\tilde{\zeta})=g(\tilde{\zeta},\tilde{\eta})=g(\tilde{\eta},\tilde{\eta})=0$, we need to show that $T(\tilde{\zeta},\tilde{\eta},\bar{\tilde{\zeta}},\bar{\tilde{\eta}})\geq0$. 

Choose $\tilde{z},\tilde{w}\in\mathbb{C}^{n+1}$ such that $\mathrm{span}_{\mathbb{C}}\{\tilde{\zeta},\tilde{\eta}\}=\mathrm{span}_{\mathbb{C}}\{\tilde{z},\tilde{w}\}$, $g(\tilde{z},\tilde{z})=g(\tilde{z},\tilde{w})=g(\tilde{w},\tilde{w})=0$ and $\tilde{z}=(z,0)$ for some $z\in\mathbb{C}^n$. By rescaling, we can normalize $\tilde{z}$ and $\tilde{w}$ such that $g(\tilde{z},\bar{\tilde{z}})=g(\tilde{w},\bar{\tilde{w}})=2$. By rotating $\tilde{w}$, we can assume $\tilde{w}=(w,ia)$ for $w\in\mathbb{C}^n$ and some nonnegative real $a$.

Hence we may write $z=e_1+ie_2$ and $w=e_3+i\lambda e_4$ for $a=\sqrt{1-\lambda^2}$. Direct calculation gives
\begin{equation*}
    T(\tilde{z},\tilde{w},\bar{\tilde{z}},\bar{\tilde{w}})=T_{1313}+\lambda^2T_{1414}+T_{2323}+\lambda^2T_{2424}-2\lambda T_{1234}+c(1-\lambda^2)\geq0.
\end{equation*}

Hence $T\in C_{PIC}$ and $\{e_1,e_2,e_3,\lambda e_4+\sqrt{1-\lambda^2}e_0\}$ is an orthonormal four-frame in $\mathbb{R}^{n+1}$ such that $T$ attains zero isotropic curvature. By Lemma 7.4 in \cite{brendle2010book}, we obtain

\begin{equation*}
    \begin{split}
        0\leq&T^{\#}(e_1,e_3,e_1,e_3)+T^{\#}(e_1,\lambda e_4+\sqrt{1-\lambda^2}e_0,e_1,\lambda e_4+\sqrt{1-\lambda^2}e_0)\\
        +&T^{\#}(e_2,e_3,e_2,e_3)+T^{\#}(e_2,\lambda e_4+\sqrt{1-\lambda^2}e_0,e_2,\lambda e_4+\sqrt{1-\lambda^2}e_0)\\
        +&2T^{\#}(e_1,e_3,\lambda e_4+\sqrt{1-\lambda^2}e_0,e_2)+2T^{\#}(e_1,\lambda e_4+\sqrt{1-\lambda^2}e_0,e_2,e_3).\\
    \end{split}
\end{equation*}

Now plug in $(T^{\#})_{ijkl}=2(\sum_{p,q=0}^{n}T_{ipkq}T_{jplq}-T_{iplq}T_{jpkq})$ and the definition of $T$ to get

\begin{equation*}
    \begin{split}
        0\leq&\sum_{p,q=0}^{n}(T_{1p1q}+T_{2p2q})\left(T_{3p3q}+\lambda^2T_{4p4q}+\lambda\sqrt{1-\lambda^2}T_{4p0q}+\lambda\sqrt{1-\lambda
        ^2}T_{op4q}+(1-\lambda^2)T_{0p0q}\right)\\
        &-\sum_{p,q=0}^n T_{12pq}\left(\lambda T_{34pq}+\sqrt{1-\lambda^2}T_{30pq}\right)\\
        &-\sum_{p,q=0}^n\left(T_{1p3q}+\lambda T_{2p4q}+\sqrt{1-\lambda^2}T_{2p0q}\right)\left(T_{3p1q}+\lambda T_{4p2q}+\sqrt{1-\lambda^2}T_{0p2q}\right)\\
        &-\sum_{p,q=0}^n\left(\lambda T_{1p4q}+\sqrt{1-\lambda^2}T_{1p0q}-T_{2p3q}\right)\left(\lambda T_{4p1q}+\sqrt{1-\lambda^2}T_{0p1q}-T_{3p2q}\right)\\
        =&\sum_{p,q=1}^n (R_{1p1q}+R_{2p2q})\left(R_{3p3q}+\lambda^2 R_{4p4q}+\frac{(1-\lambda^2)c}{2}\right)\\
        &-\sum_{p,q=1}^n \lambda R_{12pq}R_{34pq}\\
        &-\sum_{p,q=1}^n (R_{1p3q}+\lambda R_{2p4q})(R_{3p1q}+\lambda R_{4p2q})-\sum_{p,q=1}^n (\lambda R_{1p4q}-R_{2p3q})(\lambda R_{4p1q}-R_{3p2q})\\
        &+c^2(1+\lambda^2)
    \end{split}
\end{equation*}

Now since $Z\geq0$ and $Z=0$ for this particular orthonormal four-frame $\{e_1,e_2,e_3,e_4\}$ and $\lambda\in[0,1)$, we have $\frac{\partial}{\partial\lambda}Z=0$, hence $\lambda R_{1414}+\lambda R_{2424}-R_{1234}-c\lambda=0$ and $R_{1313}+R_{2323}-\lambda R_{1234}+c=0$. Consequently: 

\begin{equation*}
    \begin{split}
        &(R^2)_{1313}+\lambda^2(R^2)_{1414}+(R^2)_{2323}+\lambda^2(R^2)_{2424}-2\lambda(R^2)_{1342}-2\lambda(R^2)_{1423}\\
        =&\sum_{p.q=1}^n(R_{13pq}-\lambda R_{24pq})^2+\sum_{p,q=1}^n(\lambda R_{14pq}+R_{23pq})^2\\
        \geq&2(R_{1313}-\lambda R_{1324})^2+2(\lambda R_{2424}-R_{1324})^2
        +2(\lambda R_{1414}+R_{1423})^2+2(R_{2323}+\lambda R_{1423})^2\\
        \geq&\left[(R_{1313}-\lambda R_{1324})+(\lambda R_{2424}-R_{1324})\right]^2+\left[(\lambda R_{1414}+R_{1423})+(R_{2323}+\lambda R_{1423})\right]^2\\
        =&c^2(1+\lambda^2).
    \end{split}
\end{equation*}

We conclude that

\begin{equation*}
    \begin{split}
        &Q(R)_{1313}+\lambda^2Q(R)_{1414}+Q(R)_{2323}+\lambda^2Q(R)_{2424}-2\lambda Q(R)_{1234}\\
        =&\sum_{p,q=1}^n(R_{13pq}-\lambda R_{24pq})^2+\sum_{p,q=1}^n(\lambda R_{14pq}+R_{23pq})^2\\
        &+2\sum_{p,q=1}^n (R_{1p1q}+R_{2p2q})\left(R_{3p3q}+\lambda^2 R_{4p4q}\right)\\
        &-2\sum_{p,q=1}^n \lambda R_{12pq}R_{34pq}\\
        &-2\sum_{p,q=1}^n (R_{1p3q}+\lambda R_{2p4q})(R_{3p1q}+\lambda R_{4p2q})-2\sum_{p,q=1}^n (\lambda R_{1p4q}-R_{2p3q})(\lambda R_{4p1q}-R_{3p2q})\\
        \geq& -c(1-\lambda^2)(\Ric_{11}+\Ric_{22})
    \end{split}
\end{equation*}

$\hfill\square$

\hspace*{\fill}

We then need the following Liouville type theorem.

\begin{proposition}\label{Liouvillethm}(\cite{naber2006some})
    Let $(M^n,g,f)$ be a complete SGRS with bounded Ricci curvature. $u:M\rightarrow\mathbb{R}$ be a continuous bounded function that satisfies
    \begin{equation}\label{flapineq}
        \Delta_fu\geq0
    \end{equation}
    in the barrier sense. Then $u$ is a constant.
\end{proposition}

\textit{Proof.} Without loss of generality we assume $u\geq0$. Fix a point $x_0\in M$. For any $r>0$, one can choose a smooth cut-off function $\phi$ with support in $\{x\in M:d(x,x_0)\leq r\}$ and that $|\nabla\phi|<\frac{
C
}{r}$. Multiply $e^{-f}\phi^2u$ at both sides of equation \ref{flapineq} and use integration by parts, we get

\begin{equation*}
    \begin{split}
        0\leq&\int_Me^{-f}\phi^2u(\Delta u -\langle\nabla u,\nabla f\rangle)\\
        =&\int_M -\langle\nabla (e^{-f}\phi^2u),\nabla u\rangle-\int_Me^{-f}\phi^2u\langle\nabla u,\nabla f\rangle\\
        =&-2\int_Me^{-f}\phi u\langle\nabla\phi,\nabla u\rangle-\int_Me^{-f}\phi^2|\nabla. u|^2
    \end{split}
\end{equation*}

It follows that

\begin{equation*}
    \begin{split}
        \int_M e^{-f}\phi^2|\nabla u|^2\leq&-2\int_Me^{-f}\phi u \langle\nabla\phi,\nabla u\rangle\\
        \leq&\int_M\left(\frac{1}{2}e^{-f}\phi^2|\nabla u|^2+2e^{-f}u^2|\nabla\phi|^2\right).
    \end{split}
\end{equation*}

Hence

\begin{equation*}
    \int_M e^{-f}\phi^2|\nabla u|^2\leq 4\int_M e^{-f}u^2|\nabla\phi|^2\leq\frac{4C^2}{r^2}\int_Me^{-f}u^2.
\end{equation*}

Now let $r\rightarrow0$, since $u$ is bounded by a constant, the integral $\int_Mu^2e^{-f}$ is bounded. Therefore $|\nabla u|\equiv0$ on $M$, which means that $u$ is a constant.

$\hfill\square$

\hspace*{\fill}

\begin{proposition}\label{propPICPIC1}
    Suppose $n\geq5$, $(M,g,f)$ is a non-flat $n-$dimensional normalized SGRS that is strictly PIC. Assume further that $\nabla^2f$ is $2-$nonnegative and there exist at least one point $\bar{x}$ such that the curvature tensor at $\bar{x}$ is WPIC1. Then $(M,g)$ is WPIC1.
\end{proposition}

\textit{Proof.} Let $E$ denote the fiber bundle over $M$ such that each fiber $E_x$ is the set of all ordered orthonormal four-frame $\{e_1,e_2,e_3,e_4\}\subset T_xM$.

We define the function

\begin{equation*}
    p_1(x)=\mathrm{min}\left\{0,\inf_{\{e_1,e_2,e_3,e_4\}\in E_x,\lambda\in[0,1) }\frac{R_{1313}+\lambda^2R_{1414}+R_{2323}+\lambda^2R_{2424}-2\lambda R_{1234}}{1-\lambda^2}\right\}
\end{equation*}

 We firstly show that $p_{1}(x)>-\infty$ and is bounded below. By Proposition A.6 of \cite{Brendle19}, we know that for all $\lambda\in[0,1]$
    \begin{equation*}
        R_{1313}+\lambda^2R_{1414}+R_{2323}+\lambda^2R_{2424}-2\lambda R_{1234}+\frac{1-\lambda^2}{n-4}(\Ric_{11}+\Ric_{22}-2R_{1212})\geq0.
    \end{equation*}

Since $\nabla^2f$ is $2-$nonnegative, it follows that $0\leq\scal=\frac{n}{2}-\Delta f\leq\frac{n}{2}$.  As $R\in C_{PIC}$, there exists a constant $\Theta=\Theta(n)$ such that $|R|\leq\Theta\scal$. So

\begin{equation*}
    0\geq p_1(x)\geq-\frac{1}{n-4}|\Ric_{11}+\Ric_{22}-2R_{1212}|
\end{equation*}
is bounded below.

For fixed $x$, if $p_1(x)<0$, we claim that this infimum value must be attained by some $\{e_1,e_2,e_3,e_4\}\subset T_xM$ and some $\lambda_x\in[0,1)$. In fact, since the curvature tensor $R$ is strictly PIC, $E_x$ is compact, we know that 
\begin{equation*}
    \min_{\{e_1,e_2,e_3,e_4\}\in E_x} (R_{1313}+R_{1414}+R_{2323}+R_{2424}-2R_{1234})=\epsilon(x)>0.
\end{equation*}

Hence there exist a $\tilde{\lambda}_x$ depending on $\epsilon(x)$ and $|R|$ such that for $\lambda\in(\tilde{\lambda}_x,1)$, 
\begin{equation*}
    R_{1313}+\lambda^2R_{1414}+R_{2323}+\lambda^2R_{2424}-2\lambda R_{1234}>0.
\end{equation*}

Thus,
\begin{equation*}
    \begin{split}
        p_1(x)&=\inf_{\{e_1,e_2,e_3,e_4\}\in E_x,\lambda\in[0,1) }\frac{R_{1313}+\lambda^2R_{1414}+R_{2323}+\lambda^2R_{2424}-2\lambda R_{1234}}{1-\lambda^2}\\
        &=\min_{\{e_1,e_2,e_3,e_4\}\in E_x,\lambda\in[0,\tilde{\lambda}_x]}\frac{R_{1313}+\lambda^2R_{1414}+R_{2323}+\lambda^2R_{2424}-2\lambda R_{1234}}{1-\lambda^2}
    \end{split}
\end{equation*}
is attained by some $\{e_1,e_2,e_3,e_4\}\in E_x$ and $\lambda_x\in[0,1)$.

Now we deduce a differential inequality for $p_1(x)$. Fix $x$, suppose $p_1(x)<0$ and is attained by $\{e_1,e_2,e_3,e_4\}\in E_x$ and $\lambda_x\in[0,1)$. Choose a neighborhood $U$ of $x$ and extend $\{e_1,e_2,e_3,e_4\}$ to $U$ by parallel translation. For $y\in U$, define

\begin{equation*}
    \Bar{p}_1(y)=\frac{R_{1313}+\lambda_x^2R_{1414}+R_{2323}+\lambda_x^2R_{2424}-2\lambda_xR_{1234}}{1-\lambda_x^2}.
\end{equation*}

Then $\Bar{p}_1(y)\geq p_1(y)$ for $y\in U$ and $\Bar{p_1}(x)=p_1(x)$. Hence

\begin{equation*}
    \begin{split}
        \Delta_fp_1(x)\leq&\Delta_f\Bar{p}_1(x)=\frac{1}{1-\lambda_x^2}(\Delta_fR_{1313}+\lambda_x^2\Delta_fR_{1414}+\Delta_fR_{2323}+\lambda_x^2\Delta_fR_{2424}-2\lambda_x\Delta_fR_{1234})\\
        =&p_1(x)-\frac{1}{1-\lambda_x^2}(Q(R)_{1313}+\lambda_x^2Q(R)_{1414}+Q(R)_{2323}+\lambda^2_xQ(R)_{2424}-2\lambda_xQ(R)_{1234})
    \end{split}
\end{equation*}

Use Lemma \ref{Lem for Q(R)PIC1} for $\lambda=\lambda_x$ and $c=-p_1(x)$, we know that

\begin{equation*}
    Q(R)_{1313}+\lambda_x^2Q(R)_{1414}+Q(R)_{2323}+\lambda_x^2Q(R)_{2424}-2\lambda_xQ(R)_{1234}+(1-\lambda_x^2)p_1(x)(\Ric_{11}+\Ric_{22})\geq0
\end{equation*}

It follows that

\begin{equation}\label{diffineqforp1}
    \Delta_fp_1(x)\leq p_1(x)(1-(\Ric_{11}+\Ric_{22}))
\end{equation}
holds in the barrier sense.

By our assumption that $\nabla^2f$ is $2-$nonnegative, we know that
\begin{equation*}
    1-(\Ric_{11}+\Ric_{22})=\nabla_{11}^2f+\nabla_{22}^2f\geq0.
\end{equation*}

As a result, 
\begin{equation*}
    \Delta_fp_1(x)\leq0.
\end{equation*}

Since $p_1(x)$ is bounded, and $0\leq|Ric|\leq C(n)\scal=C(n)(\frac{n}{2}-\Delta f)\leq \frac{C(n)n}{2}$ is bounded by a constant. Use Proposition \ref{Liouvillethm} for $u=-p_1$, we know that $p_1(x)$ is a constant.

Now since there exist at least one point $\bar{x}$ on $M$ such that the curvature tensor is WPIC1 at $\bar{x}$, we know that $p_{1}\equiv 0$. This means that $(M,g)$ is WPIC1.

$\hfill\square$

\hspace*{\fill}

Now we are ready to proof Theorem \ref{main thm}. By Lemma 4.2 of \cite{bamler2019} or Proposition 6.2 of \cite{Li2019}, we conclude that $(M,g)$ is WPIC2. Notice that the Ricci curvature tensor is nonnegative, and the associated Ricci flow of $(M,g,f)$ satisfies the evolution equation

\begin{equation*}
    \frac{\partial\Ric_{ij}}{\partial t}=\Delta\Ric_{ij}+2\sum_{k,l}R_{ikjl}\Ric_{kl}.
\end{equation*}

The term $(R\star \Ric)_{ij}=\sum_{k,l}R_{ikjl}\Ric_{kl}$ is nonnegative definite. By strong maximum principle, the null space of $\Ric$ is invariant under parallel translation. Moreover, the null space of $\Ric$ must have $\mathrm{rank}=0$ or $1$, otherwise this contradicts with the condition that $(M,g)$ is PIC and WPIC2.

If the null space of $\Ric$ has rank 0, then $\Ric$ is strictly positive. Theorem 2 of \cite{MW15} shows that $M$ is compact. By Proposition 6.6 of \cite{Brendle19}, $(M,g)$ is strictly PIC2, and hence $(M,g)$ is a quotient of $S^{n}$ by \cite{BrendleSchoen07}. 

If the null space of $\Ric$ has rank 1, then the universal cover $\tilde{M}$ of $M$ splits as product of an $(n-1)-$dimensional SGRS $(N,g_N,f_N)$ with $\mathbb{R}$. Moreover, $(N,g_N)$ has strictly positive Ricci. By the argument above, $(N,g_N)$ is an $(n-1)-$dimensional spherical space form. Hence $(M,g)$ is a quotient of $S^{n-1}\times\mathbb{R}$. This completes the proof.

\bibliographystyle{alpha}
\bibliography{PICSGRSbib}

\end{document}